\documentclass[11pt,twoside]{article}

\usepackage{amsfonts}
\usepackage{amsmath}
\usepackage{amssymb}
\usepackage{amsthm}
\usepackage{graphics} 
\usepackage{epsfig} 
\usepackage{latexsym}

\graphicspath{{./figures/}} 

\usepackage[textwidth=16cm,textheight=22cm,top=3cm,bottom=3cm]{geometry}
\usepackage{fancyhdr}

\newtheorem{theorem}{Theorem}
\newtheorem{lemma}{Lemma}
\newtheorem{corollary}{Corollary}

\newtheorem{remark}{Remark}
\theoremstyle{definition}

\usepackage{color}

\newcommand{\co}{\circledast}

\newcommand{\chil}{\raisebox{0.5mm}{\large $\chi$}}
\newcommand{\qtextq}[1]{\quad\text{#1}\quad}

\newcommand{\drop}[1]{}
\newcommand{\no}{\noindent}
\newcommand{\fer}[1]{(\ref{#1})}
\newcommand{\qtext}[1]{\quad\text{#1}}

\newcommand{\eps}{\varepsilon}
\newcommand{\vfi}{\varphi}

\newcommand{\grad}{\nabla}
\renewcommand{\o}{\omega}
\newcommand{\p}{\partial}

\newcommand{\N}{\mathbb{N}}
\newcommand{\R}{\mathbb{R}}

\def\O{\Omega}

\newcommand{\abs}[1]{| #1 |}
\newcommand{\norm}[1]{\| #1 \|}

\title{Symmetrization in nonlocal diffusion problems   
\thanks{Supported by Spanish MCI Project MCI-21-PID2020-116287GB-I00. 
}}

\author{Gonzalo Galiano\thanks{Department of Mathematics, University of Oviedo, Spain 
  ({\tt galiano@uniovi.es}).}
}

\date{}

\begin{document}

\maketitle

\begin{abstract}
We study Talenti's type symmetrization properties for solutions of linear stationary and evolution problems.

Our main result establishes the comparison in norm between the solution of  a problem and its symmetric version when nonlocal diffusion defined through integrable kernels is replacing the usual local diffusion defined by a second order differential operator. 

Using an approximation argument, we recover, as a corollary of our results, the classical Talenti's theorem. A novelty of our approach is that we replace  the measure geometric tools employed in Talenti's  proof by the use of the Riesz's rearrangement inequality, giving thus an alternative and somehow simpler proof than Talenti's one.\vspace{0.25cm}

\no\emph{Keywords: }Nonlocal diffusion, Schwarz's symmetrization, Talenti's theorem, {Riesz's inequality.}

\end{abstract}


\section{Introduction}

In this introduction, we briefly describe our results for the stationary problem. As we shall see, similar results are deduced for the evolution  problem.

Let $\O\subset\R^N$ be an open and bounded set, and  let $f\in L^p(\O)$, with $p\geq1$, be non-negative. The \emph{kernel} $J:\R^N\to\R$ is assumed to be a non-negative radially non-increasing function with $\norm{J}_{L^1(\R^N)}=1$. 

Let $v\in L^p(\R^N)$ be a non-negative solution of the nonlocal stationary problem with homogeneous Dirichlet boundary condition, this is, satisfying
\begin{align*}
&-\int_{\R^N} J(x-y) (v(y)-v(x))dy +c v(x)=  f(x)&& \text{for }x\in \O,& \\
& v(x) = 0 &&\text{for } x\in\R^N\backslash \O,&
\end{align*}
where $c\geq0$ is a constant. 

We consider the problem of establishing a comparison property between $v$ and the solution of the corresponding \emph{symmetrized problem}
\begin{align*}
&-\int_{\R^N} J(x-y) (u(y)-u(x))dy +cu(x)=  f^*(x)&& \text{for }x\in \O^*,& \\
& u(x) = 0 &&\text{for } x\in\R^N\backslash \O^*.&
\end{align*}
Here, $\O^*$ denotes the ball of $\R^N$ centered at zero and with the same volume as $\O$. The function $f^*$ is the Schwarz's symmetrization of $f$, which is a radially non-increasing function having the same level sets measures as $f$.

Our first result establishes the comparison in norm,
\begin{align*}
\norm{v}_{L^p(\O)}	\leq \norm{u}_{L^p(\O^*)}.
\end{align*}
This property, satisfied in the nonlocal diffusion case, may be extended to the local diffusion case by an approximation argument. Consider the rescaled kernel $J_\eps(x) = C_1\eps^{-(N+2)}J(x/\eps)$, for $\eps>0$  and $C_1$ given in \fer{def:C1}. Then, under suitable regularity assumptions on $\O$ and $f$, we have 
\begin{align*}
\norm{v_\eps - V}_{L^p(\O)}\to 0 \qtextq{as}\eps\to0,	
\end{align*}
where $v_\eps$ is the solution of the nonlocal stationary problem with kernel $J_\eps$ and $V$ is the solution of 
\begin{align*}
&-\Delta V + cV=  f&& \text{in } \O,& \\
& V = 0 &&\text{on } \p\O.&
\end{align*}
We may use the same rescaled kernel to define a sequence, $u_\eps$, for the symmetrized nonlocal stationary problem, enjoying the property $\norm{u_\eps - U}_{L^p(\O)}\to 0$ as $\eps\to0$, where $U$ is the solution of the symmetrized  problem
\begin{align*}
&-\Delta U + cU=  f^*&& \text{in } \O^*,& \\
& U = 0 &&\text{on } \p\O^*.&
\end{align*}

Our second result is then clear. Since 	$\norm{v_\eps}_{L^p(\O)}	\leq \norm{u_\eps}_{L^p(\O^*)}$, $v_\eps\to V$ in $L^p(\O)$ and 
$u_\eps\to U$ in $L^p(\O^*)$, we deduce that 
\begin{align*}
\norm{V}_{L^p(\O)} \leq 	\norm{U}_{L^p(\O^*)},
\end{align*}
which is a weak form of Talenti's theorem (Talenti \cite{Talenti1976} provides the stronger pointwise estimate $V^*(x) \leq U(x)$, for $x\in\O^*$).

The article is organized as follows.  In Section~\ref{main}, we introduce some common tools related to Schwarz's symmetrization and state our assumptions and main results. In Sections~\ref{sec:stationary} and \ref{sec:evolution}, we prove the results related to the stationary and the evolution problems, respectively.

\section{Main results}\label{main}


\subsection{Schwarz's symmetrization}

Let $E\subset\R^N$ be a measurable set of finite measure, and let $\chil_E:\R^N\to\R$ be its characteristic  function, i.e., defined by $\chil_E(x)=1$ if $x\in E$ and  $\chil_E(x)=0$ otherwise. 

The \emph{symmetric rearrangement of $E$} is 
the ball $E^*\subset\R^N$ centered at zero with  $\abs{E^*}=\abs{E}$, i.e., with radius $(\abs{E}/\o_N)^{1/N}$, where $\o_N$ denotes the volume of the $N-$dimensional unit ball.

For a non-negative measurable function $f:\R^N\to\R$ vanishing at infinity, the \emph{Schwarz's symmetrization of $f$} is
\begin{align*}
f^*(x) = \int_0^\infty \chil_{\{f>s\}^*}(x)ds,	
\end{align*}
where, by definition, $(\chil_E)^* = \chil_{E^*}$. Thus, the level sets of $f^*$ are the rearrangements of the level sets of $f$, implying the equimeasurability property
\begin{align*}
\abs{\{x:f^*(x)>s\}}	= \abs{\{x:f(x)>s\}}.
\end{align*}

The Schwarz's symmetrization of a function inherits many measure geometric properties from its source function $f$. It also fulfils some optimization properties with respect to integration. In this article we shall make a recurrent use of Hardy-Littlewood's inequality
\begin{align*}
\int_{\R^N} f_1(x)f_2(x)dx \leq 	\int_{\R^N} f_1^*(x)f_2^*(x)dx,
\end{align*}
and of  Riesz's inequality \cite{Riesz1930}
\begin{align*}
\int_{\R^N} f_1(x) \Big(\int_{\R^N} f_2(x-y)f_3(y)dy \Big) dx \leq 	
\int_{\R^N} f_1^*(x) \Big(\int_{\R^N} f_2^*(x-y)f_3^*(y)dy \Big) dx,
\end{align*}
where $f_1,f_2,f_3$  are measurable non-negative functions vanishing at infinity. We shall also use the generalized Riesz's inequality deduced in \cite{Brascamp1974, Lieb2001}.

There is an extensive literature on Schwarz's symmetrization, as well as on  other types of symmetrizations and functional rearrangements. We refer the reader to the books \cite{Polya1951, Bandle1976, Kesavan2006} for a detailed description of these tools, and to the review article \cite{Talenti2016} as an abundant source of applications and references.


\subsection{The problems}

Let $E\subset\R^N$ be an open and bounded set with smooth boundary, let $h:E\to\R$ be a non-negative measurable function, and let $\eps>0$. The \emph{stationary nonlocal diffusion problem}, SNP($E,h,\eps$), consists on finding $w:\R^N\to\R$ such that
\begin{align*}
&- \int_{\R^N} J_\eps(x-y) (w(y)-w(x))dy +cw(x)=  h(x)&& \text{for }x\in E,& \\
& w(x) = 0 &&\text{for }x\in\R^N\backslash E.& 
\end{align*}
The corresponding \emph{stationary local diffusion problem}, SLP($E,h$), consists on finding $W:E\to\R$ such that 
\begin{align*}
&- \Delta W +cW=  h&& \text{in } E,&\\
& W = 0 &&\text{on } \p E.&
\end{align*}

Accordingly, the \emph{evolution nonlocal diffusion problem}, ENP($E,h,\zeta,\eps$), with $E$ and $\eps$ like above, $h:(0,T)\times E\to\R$ and  $\zeta:E\to\R$  non-negative measurable functions, consists on finding $w:(0,T)\times \R^N\to\R$ such that
\begin{align*}
&\p_t w(t,x)-  \int_{\R^N} J_\eps(x-y) (w(t,y)-w(t,x))dy +cw(t,x)=  h(t,x)&& & 
\end{align*}
for $(t,x)\in (0,T)\times E$, and 
\begin{align*}
& w(t,x) = 0 &&\text{for }(t,x)\in (0,T)\times\R^N\backslash E,& \\
& w(0,x) = \zeta(x) &&\text{for }x\in E.& 
\end{align*}
The corresponding \emph{evolution local diffusion problem}, ELP($E,h,\zeta$), consists on finding $W:(0,T)\times E\to\R$ such that 
\begin{align*}
&\p_tW - \Delta W +cW=  h && \text{in } (0,T)\times E,&\\
& W = 0 &&\text{on } (0,T)\times \p E,&\\
& W(0,\cdot) = \zeta &&\text{in } E.&
\end{align*}



\subsection{Assumptions}

The minimal hypothesis for proving the well-posedness of the nonlocal diffusion and the local diffusion problems are different. Proving the convergence of the solutions of the nonlocal diffusion rescaled problems to their local diffusion versions requires  additional assumptions. 

Since in this article we are interested in the comparison results between solutions of a problem and of its symmetrized version, we give here just the assumptions needed to prove them. 

The additional hypotheses needed to ensure the existence, uniqueness and regularity of solutions, as well as the convergence of solutions of nonlocal diffusion rescaled problems to their local diffusion counterpart, are assumed to hold.  
The reader is referred to \cite{Andreu2010,Evans1998} for details on these questions.

\medskip

\no\textbf{Assumptions (H)}
\begin{enumerate}
\item $ J \in L^1(\R^N)\cap L^r(\R^N)$, for some $1<r\leq 2$, is a non-negative function which is \emph{radially non-increasing}, this is, such that  (i) $J(x)=J(y)$ if $\abs{x}=\abs{y}$,  and (ii) $J(x)\leq J(y)$ if $\abs{x}\geq\abs{y}$. Besides, we assume $\norm{J}_{L^1(\R^N)}=1$ and define the constant 
\begin{align}
\label{def:C1}
C_1 = 2/\sigma^2 \qtextq{with}\sigma^2 = \int_{\R^N} \abs{x}^2 J(x) dx <\infty.
\end{align}
\item The data, $\O\subset\R^N$, $F:\O\to\R$, $G:(0,T)\times\O\to\R$, and $v_0:\O\to\R$ are regular enough to provide unique solutions with $L^p$ regularity, for some $1\leq p\leq\infty$, to problems SNP($\O,F,\eps$), SLP($\O,F$), ENP($\O,G,v_0,\eps$), SNP($\O,G,v_0$), and their corresponding symmetrized versions. We always assume that, at least, $F\in L^p(\O)$ and $G\in L^p((0,T)\times\O)$.
\item In addition to the (H)$_2$, we assume that the the regularity of the data is  enough to provide the convergence of solutions of nonlocal diffusion rescaled problems to their local diffusion counterpart. Explicitly, if $v_\eps\in L^p(\R^N)$ and $u_\eps\in L^p(\R^N)$ are the solutions of  SNP($\O,F,\eps$) and SNP($\O^*,F^*,\eps$), and if $V\in L^p(\O)$ and $U\in L^p(\O^*)$ are the solutions of  SLP($\O,F,\eps$) and SLP($\O^*,F^*,\eps$) then we have
\begin{align*}
	\norm{v_\eps-V}_{L^p(\O)}\to0 \qtextq{and} \norm{u_\eps-U}_{L^p(\O^*)}\to0.
\end{align*}
A similar property is assumed for the evolution problems.
\item This assumption is related to the convergence of time discrete approximations schemes to their corresponding time continuous versions. See \fer{h4.1}-\fer{h4.2} for details.
\end{enumerate}

\begin{remark}
The condition $ J \in L^r(\R^N)$ for $1<r\leq 2$ is unusual, although not specially limitating. The use of the central limit theorem in the proof of Theorem~\ref{th:stationary} requires it.
\end{remark}


\subsection{Main results}

We obtain similar comparison results for the stationary and the evolution problems. For the stationary problems, we have:
\begin{theorem}
\label{th:stationary}
Assume (H) and set $\eps>0$.  Let $v_\eps \in L^p(\R^N)$ and  $u_\eps \in L^p(\R^N)$ be the solutions of SNP($\O,F,\eps$) and SNP($\O^*,F^*,\eps$), respectively. Then
$u_\eps$ is radially non-increasing and  
\begin{align}
\label{th:comp}
\norm{v_\eps}_{L^p(\O)}	\leq \norm{u_\eps}_{L^p(\O^*)}.
\end{align}
\end{theorem}

\begin{corollary}
\label{cor:stationary}
Let $V$ be the solution of SLP($\O,F$) and  let $U$ be the solution of SLP($\O^*,F^*$). Then 
$U$ is radially non-increasing and  
\begin{align*}
\norm{V}_{L^p(\O)}	\leq \norm{U}_{L^p(\O^*)}.
\end{align*}
\end{corollary}

For the evolution problems, we have:
\begin{theorem}
\label{th:evolution}
Assume (H) and set $\eps>0$. Let $v_\eps \in L^p((0,T)\times\R^N)$ and  $u_\eps \in L^p((0,T)\times\R^N)$ be the solutions of ENP($\O,G,v_0,\eps$) and ENP($\O^*,G^*,v_0^*,\eps$), respectively. Then $u_\eps(t,\cdot)$ is radially non increasing for a.e. $t\in(0,T)$ and 
\begin{align}
\label{th2:comp}
\norm{v_\eps}_{L^p((0,T)\times\O)}	\leq \norm{u_\eps}_{L^p((0,T)\times\O^*)}.
\end{align}
\end{theorem}

\begin{corollary}
\label{cor:evolution}
Let $V$ be the solution of ELP($\O,G,v_0$) and  let $U$ be the solution of SLP($\O^*,G^*, v_0^*$). Then 
$U(t,\cdot)$ is radially non increasing for a.e. $t\in(0,T)$ and 
\begin{align*}
\norm{V}_{L^p((0,T)\times\O)}	\leq \norm{U}_{L^p((0,T)\times\O^*)}.
\end{align*}
\end{corollary}

\begin{remark}
\label{rem:cor}
Under assumption (H)$_3$, the corollaries are straightforward to prove. For instance, in the stationary case, 
using the triangle inequality and Theorem~\ref{th:stationary} we obtain
\begin{align*}
	\norm{V}_{L^p(\O)} &\leq \norm{V-v_\eps}_{L^p(\O)} + 	\norm{v_\eps}_{L^p(\O)} \leq 	
	 \norm{V-v_\eps}_{L^p(\O)} + 	\norm{u_\eps}_{L^p(\O^*)} \\ 
	 & \leq  \norm{V-v_\eps}_{L^p(\O)} + 	\norm{u_\eps-U}_{L^p(\O^*)} + 	\norm{U}_{L^p(\O^*)}, 
\end{align*}
and letting $\eps\to0$, we deduce the assertion. The evolution case is treated similarly.
\end{remark}

\subsection{Notation }

 

For a mesurable set $E\subset\R^N$ and a function $\vfi:E\to\R$, we define  the \emph{extension by zero of $\vfi$ to $\R^N$} as 
\begin{align*}
\overline{\vfi}(x) = \begin{cases}
 \vfi(x)& \text{if }x\in E,\\
 0 & \text{if }x\in \R^N \backslash  E.	
 \end{cases}
\end{align*}
For $\vfi,\psi\in L^1(\R^N)$, we denote the convolution of $\vfi$ with $\psi$ in $\R^N$ by the usual symbol $*$, this is, 
\begin{align*}
\vfi*\psi (x) = \int_{\R^N} \vfi(x-y) \psi(y)dy.
\end{align*}
We recall that the convolution in $\R^N$ is conmmutative. 
 For the convolution of $\vfi$ with $\psi$ in $E$, we use the symbol $\co_E$
\begin{align*}
\vfi\co_E\psi (x) = \int_{E} \vfi(x-y) \psi(y)dy.
\end{align*}
If the context is clear, we just write $\co$.
Observe that if $\vfi,\psi\in L^1(\R^N)$ with $\psi=0$ in $\R^N \backslash  E$ then $\vfi*\psi = \vfi \co_E\psi$. 


We write $(J*)^1 = J $, and denote by $(J*)^k$, for $k=1,2,\ldots$, to the recurrent convolution in $\R^N$ of $k$ kernels, $J$. This is, $(J*)^2=J*J$, $(J*)^3=J*J*J$, and so on.  Observe that $\norm{(J*)^k}_{L^1(\R^N)}=\norm{J}_{L^1(\R^N)}^k$. Correspondingly, we write $(J\co_E)^k$ to denote  the convolution in $E$ of $k$ kernels, $J$. In this case,  $\norm{(J\co_E)^k}_{L^1(\O)}\leq \norm{J}_{L^1(\O)}^k$.


\section{Proof of Theorem~\ref{th:stationary}}\label{sec:stationary}

Problem SNP($E,h,\eps$) may be reformulated as 
\begin{align}
&w(x)= \alpha \int_{\R^N} \rho_\epsilon (x-y) w(y)dy +  \xi(x)&& \text{for }x\in E , & \label{jope1} \\
& w(x) = 0 &&\text{for }x\in\R^N\backslash E,& \label{jope2}
\end{align}
where 
$\alpha=1/(1+c\eps^2)$, $\xi=\eps^2 h/(1+c\eps^2)$, and $\rho_\eps  = \eps^2 J_\eps$. Observe that $\norm{\rho_\eps }_{L^1(\R^N)} =  \norm{J}_{L^1(\R^N)}=1$, for all $\eps>0$.

\medskip

The limit $\eps\to0$ is only considered for proving Corollaries~\ref{cor:stationary} and \ref{cor:evolution}. Since in the proofs of Theorems~\ref{th:stationary} and \ref{th:evolution} the parameter $\eps$ is kept fixed, we drop it from the subindices of the different functions depending on it in order to have a cleaner notation.

We refer to  problem \fer{jope1}-\fer{jope2} as to  AUX($E,\xi$).
Let $v$ and $u$ be the solutions of AUX($\O,f$) and AUX($\O^*,f^*$), respectively, for  $f= \eps^2 F/(1+c\eps^2)$. Thus, $v$ and $u$ are the solutions of SNP($\O,F,\eps$) and SNP($\O^*,F^*,\eps$), respectively.

The proof of the theorem uses two lemmas. We first prove the main result stating the lemmas when  required and, afterwards, we prove the lemmas. We start with the comparison property stated in the theorem.

Let $\vfi_0\in L^{p'}(\O)$ be a non-negative function. Multiplying the first equation of AUX($\O,f$)  by $ \vfi_0$, integrating in $\O$ and using the symmetry of $\rho  $ to interchange the convolved functions, we obtain
\begin{align}
\int_\O v(x)\vfi_0(x)dx  
&=	\alpha\int_\O v(y) \int_{\O}\rho  (x-y) \vfi_0(x) dx dy + \int_\O f(x)\vfi_0(x)dx \nonumber \\
&=	\int_\O v(x) \vfi_1(x) dx + \int_\O f(x)\vfi_0(x)dx, \label{p1:2}
\end{align}
with $\vfi_1 = \alpha \rho  \co \vfi_0$. 
Multiplying now the first equation of AUX($\O,f$)  by $\vfi_1$ we, similarly, get
\begin{align*}
\int_\O v(x)\vfi_1(x)dx =  \int_\O v(x) \vfi_2(x) dx +  \int_\O f(x)\vfi_1(x)dx,
\end{align*}
with $\vfi_2 = \alpha^2(\rho  \co)^2 \co  \vfi_0$. Thus,  replacing in \fer{p1:2}, we obtain
\begin{align*}
\int_\O v(x)\vfi_0(x)dx = \int_\O v(x) \vfi_2(x) dx +  \int_\O f(x)(\vfi_0(x)+\vfi_1(x))dx.
\end{align*}
Therefore, the following identity holds for $k=1,2,\ldots$
\begin{align}
\label{p1:3}
\int_\O v(x)\vfi_0(x)dx = \int_\O v(x) \vfi_{k}(x) dx +  \sum_{i=0}^{k-1}\int_\O f(x)\vfi_i(x)dx,
\end{align}
with $\vfi_i = \alpha^i(\rho  \co)^i \co \vfi_0$.

In a similar way, testing the first equation of AUX($\O^*,f^*$)   with $\vfi_0^*$, we obtain 
\begin{align}
\label{p1:4}
\int_{\O^*} u(x)\vfi_0^*(x)dx = \int_{\O^*} u(x) \psi_{k}(x) dx +  \sum_{i=0}^{k-1}\int_{\O^*} f^*(x)\psi_i(x)dx,
\end{align}
with $\psi_i = \alpha^i(\rho  \co)^i \co  \vfi_0^*$. Notice that Lemma~\ref{lemma:convolution} ensures that $\psi_i$ is radially non-increasing.
Using  Riesz's inequality, we deduce 
\begin{align*}
\int_\O f(x)\vfi_i(x)dx &= \alpha^i \int_{\O}  f(x)\int_{\O} (\rho  \co)^i (x-y)  \vfi_0(y)dy dx	\\
&\leq \alpha^i \int_{\O^*} f^*(x)\int_{\O^*} (\rho  \co)^i (x-y)   \vfi_0^*(y)dy dx
= \int_{\O^*} f^*(x)\psi_i(x)dx.
\end{align*}
This estimate in combination with \fer{p1:3} and \fer{p1:4} yields
\begin{align}
\label{p1:5}
\int_\O v(x)\vfi_0(x)dx \leq \int_{\O^*} u(x)\vfi_0^*(x)dx+ \int_\O v(x) \vfi_{k}(x) dx - \int_{\O^*} u(x) \psi_{k}(x) dx.
\end{align}
For dealing with the last two terms of this expression, we use the following lemma.
\begin{lemma}
\label{lemma:CLT}
Let $E\subset\R^N$ be a bounded measurable set, $\psi\in L^p(E)$ and $\vfi\in L^{p'}(E)$, for $1\leq p\leq\infty$. Then 
\begin{align}
\label{lemma:clt1}
	\lim_{i\to\infty}\int_E \vfi(x) (\rho  *)^i * \bar \psi(x) dx =0.
\end{align}	
The same convergence is attained when replacing $*$ by $\co_E$.
\end{lemma}
Therefore, using the lemma in \fer{p1:5} we obtain,  in the limit $k\to\infty$,
\begin{align*}
\int_\O v(x)\vfi_0(x)dx \leq \int_{\O^*} u(x)\vfi_0^*(x)dx.
\end{align*}
For $p=1$, we choose $\vfi_0=1$ to obtain $\norm{v}_{L^1(\O)}\leq \norm{u}_{L^1(\O^*)}$.
For $p>1$, recalling that 
$(\Phi(v))^*=\Phi(v^*)$ for any non-decreasing function $\Phi$
we deduce, for $\Phi(s)=s^{p-1}$,
\begin{align}
\label{p1:7}
\int_{\O^*} \abs{v(x)}^pdx \leq \int_{\O^*} u(x)\abs{v^*(x)}^{p-1}dx.
\end{align}
Finally,  we deduce \fer{th:comp} by  using H\"older's inequality and the elementary  property
$\norm{v}_{L^p(\O)} = \norm{v^*}_{L^p(\O^*)}$.

\bigskip

We continue the proof of the theorem showing that the solution of the symmetrized problem AUX($\O^*,f^*$), and thus of the stationary problem SNP($\O^*,F^*,\eps$), is radially non-increasing. We shall use the following elementary result.
\begin{lemma}
	Let $B_R\subset\R^N$ be the ball of radius $R>0$ centered at the origin. Let  $\vfi,\psi\in L^1(B_R)$ be radially non-increasing functions. Then
$\vfi\co_{B_R}\psi$ is radially non-increasing. The same remains true replacing $B_R$ by $\R^N$ and $\co_{B_R}$ by $*$. 
\label{lemma:convolution}
\end{lemma}

Let $u_0\in L^1(\O^*)$ be a radially non-increasing function, and consider the sequence $u_k:\O^*\to\R$ given by, for $k=1,2,\ldots$, and $x\in\O^*$,
\begin{align}
\label{sol.implicita4}
	u_{k}(x)= \alpha \rho  \co_{\O^*}  u_{k-1}(x) + f^*(x) .
\end{align}
By Lemma~\ref{lemma:convolution}, $u_{k}$ is radially non-increasing. 
Solving the recursivity  \fer{sol.implicita4}, we obtain  
\begin{align}
\label{sol.implicita3}
 u_{k}(x) = \alpha^k (\rho  \co)^{k} \co u_0(x)+  \Big( f^*(x) + \sum_{i=1}^{k-1} \alpha^i(\rho  \co)^i * f^*(x) \Big),\qtext{for }k=1,2,\ldots
\end{align}
For $x\in\O^*$, the solution of AUX($\O^*,f^*$) satisfies
\begin{align}
\label{sol.implicita2}
u(x) = \alpha \rho  \co u(x)+ f^*(x) ,
\end{align}
where we used that $u=0$ in $\R^\N\backslash\O^*$. 
Therefore, for $k=2,3,\ldots$ and $x\in\O^*$, we deduce the implicit formula 
\begin{align}
u(x) = \alpha^k (\rho  \co)^k \co u(x)+  \Big( f^*(x) +  \sum_{i=1}^{k-1} \alpha^i(\rho  \co)^i \co f^*(x)\Big).	
\label{sol.implicita}
\end{align}
From \fer{sol.implicita3} and \fer{sol.implicita}, we obtain, for $x\in \O^*$ and $k=2,3,\ldots$, 
\begin{align*}
 u_{k}(x) = u(x)+ \alpha^k(\rho  \co)^{k} \co( u_0-u) (x).
\end{align*}
Defining $ u_0 =  f^*$ and using \fer{sol.implicita2} we obtain, for $x\in \O^*$ and $k=2,3,\ldots$ 
\begin{align*}
 u_{k}(x) = u(x) - \alpha^{k+1}(\rho  \co)^{k+1} \co u(x),
\end{align*}
implying 
\begin{align}
\label{rad1}
	\norm{u_{k}- u}_{L^1(\O^*)} = \alpha^{k+1}\int_{\O^*} (\rho  \co )^{k+1} \co u (x)dx.
\end{align}


Using Lemma~\ref{lemma:CLT} with $E=\O^*$, $\vfi =1$ and $\psi = u|_{\O^*}$ in \fer{rad1}
we deduce that
$u_k\to u|_{\O^*}$ strongly in $L^1(\O^*)$. At least for a subsequence, we have that $u_{k_j}\to u|_{\O^*}$ a.e. in $\O^*$. Since $u_{k_j}$ are radially non-increasing, it follows  that the solution $u$ of AUX($\O^*,f^*$) is also radially non-increasing. This finishes the proof of Theorem~\ref{th:stationary}. $\hfill\Box$

\medskip

\no\textbf{Proof of Lemma~\ref{lemma:CLT}. }
We have
\begin{align*}
	\left|\int_E \vfi(x) (\rho  *)^k * \bar \psi(x) dx \right|\leq \norm{\vfi }_{L^{p'}(E)}\norm{(\rho  *)^k * \bar \psi }_{L^p(E)}.
\end{align*}
Since  $E$ is bounded, there exists a ball $B_R$ such that $E\subset B_R$ and $x-y \in B_R$ for all $x,y\in E$. Like in the case $E=\R^N$, we have $\norm{(\rho  *)^k * \bar \psi }_{L^p(E)}\leq \norm{(\rho  *)^k}_{L^1(B_R)}\norm{\psi }_{L^p(E)}$. Indeed,
\begin{align*}
\int_E &\left|\int_{\R^N} (\rho  *)^k (x-y)\bar\psi(y)dy \right|^p dx  = 	
\int_E \left|\int_{E} (\rho  *)^k(x-y)\psi(y)dy \right|^p dx  \\
& =
\int_{\R^N} \left|\int_{\R^N} (\rho  *)^k(x-y)\bar\psi(y) \chil_E(y)\chil_{B_R}(x-y)dy\right|^p \chil_E(x)dx \nonumber \\
&=  \norm{ \big((\rho  *)^k\chil_{B_R}\big) * (\psi \chil_E)}_{L^p(\R^N)}^p  \\
& \leq 
\norm{ (\rho  *)^k\chil_{B_R}}_{L^1(\R^N)}^p \norm{\psi \chil_E}_{L^p(\R^N)}^p . 
\end{align*}


According to the central limit theorem  \cite[\textsection46, Th.~1]{Gnedenko1968}, we have 
\begin{align}
\label{th:cl}
	\sigma\sqrt{k} (\rho  *)^{k}(\sigma\sqrt{k}x) \to \Phi(x) \qtext{as }k\to \infty\qtext{uniformly for }x \text{ in }\R^N, 
\end{align}
where $\Phi$ is the Gaussian of zero mean and identity covariance matrix. Let, for $x\in\R^N$, 
\begin{align*}
U_k(x) = \sigma\sqrt{k}(\rho  *)^{k}(\sigma\sqrt{k}x) \chil_{B_{r_k}}(x) dx ,
\end{align*}
with $r_k = R/(\sigma\sqrt{k})$, so that $\int_{B_{R}} (\rho  *)^{k} = \int_{\R^N}  U_k$. The uniform convergence \fer{th:cl} implies that for all $\delta>0$ there exists $K\in\N$ such that $\abs{\sigma\sqrt{k}(\rho  *)^k (0) - \Phi(0)}<\delta$ for $k>K$. Since $(\rho  *)^k$ is radially non-increasing, we deduce
\begin{align*}
U_k(x) \leq (\delta+\Phi(0))	\chil_{B_{r_k}}(x) \leq (\delta+\Phi(0))	\chil_{B_{r_1}}(x).
\end{align*}
The first inequality shows that $U_k(x)\to 0$ for a.e. $x\in\R^N$, while the second implies  that $U_k$ is dominated by $(+\Phi(0))	\chil_{B_{r_1}}\in L^1(\R^N)$. Therefore, the theorem of dominated convergence ensures that $\norm{(\rho  *)^k}_{L^1(B_R)}\to0$ as $k\to\infty$, and \fer{lemma:clt1} follows. 

Finally, on noting that 
\begin{align*}
\norm{(\rho  \co )^k \co \bar \psi }^p_{L^p(E)} & = 	
\int_E \left|\int_{E} (\rho  \co)^k(x-y)\psi(y)dy \right|^p dx  \\
& \leq 
\norm{ (\rho  \co)^k}_{L^1(B_R)}^p \norm{\psi}_{L^p(E)}^p 
\leq 
\norm{ (\rho  *)^k}_{L^1(B_R)}^p \norm{\psi}_{L^p(E)}^p , 
\end{align*}
we see that the same arguments may be employed when replacing $*$ by $\co$.
$\hfill\Box$

\medskip

\no\textbf{Proof of Lemma~\ref{lemma:convolution}.} We first prove it for the convolution of functions defined in $\R^N$. The result in $B_R$ follows by taking the restriction to $B_R$ of the extensions $\bar\vfi,\bar\psi$.

We start by checking that $A(x)=\vfi*\psi(x)$ is radial.
Let  $x_1,x_2\in \R^N$ be such that $\abs{x_1}=\abs{x_2}$. Then, there exists a orthogonal rotation matrix, $G$, such that 
$x_2=Gx_1$ and $\abs{\det G}=1$. Introducing the change of integration variable $y=Gz$ and taking into account that  $\abs{Gz}=\abs{z}$, we obtain 
\begin{align*}
A(x_2)=\int_{\R^N} \vfi(G(x_1-z)) \psi(Gz)dz	= \int_{\R^N} \vfi(x_1-z)\psi(z)dz =A(x_1),
\end{align*}
so that $A$ is radial.

We now check that $A$ is radially non-increasing. We start assuming that $\psi\in W^{1,1}(\R^N)$.  Then 
\begin{align*}
 \p_{x_j}A(x) &= \int_{\R^N}  \bar  \vfi_R(\abs{y})  \bar \psi_R'(\abs{x-y})\frac{x_j-y_j}{\abs{x-y}}dy = \int_{\R^N}  \bar  \vfi_R(\abs{x-z}) \bar  \psi_R'(\abs{z})\frac{z_j}{\abs{z}}dz,
\end{align*}
and thus 
\begin{align*}
\grad A(x) \cdot x &= \int_{\R^N}  \bar  \vfi_R(\abs{x-z}) \bar  \psi_R'(\abs{z})\frac{z\cdot x}{\abs{z}}dz \\
&= \int_{z\cdot x >0}  \bar  \vfi_R(\abs{x-z}) \bar  \psi_R'(\abs{z})\frac{z\cdot x}{\abs{z}}dz + \int_{z\cdot x <0}  \bar  \vfi_R(\abs{x-z}) \bar  \psi_R'(\abs{z})\frac{z\cdot x}{\abs{z}}dz \\
&= \int_{z\cdot x >0}  \bar  \vfi_R(\abs{x-z}) \bar  \psi_R'(\abs{z})\frac{z\cdot x}{\abs{z}}dz - \int_{z\cdot x >0}  \bar  h(\abs{x+z}) \bar  \psi_R'(\abs{z})\frac{z\cdot x}{\abs{z}}dz \\
&= \int_{z\cdot x >0}  \big(\bar  \vfi_R(\abs{x-z})-\bar  \vfi_R(\abs{x+z})\big) \bar  \psi_R'(\abs{z})\frac{z\cdot x}{\abs{z}}dz  .
\end{align*}
For $z\in\R^N$ we have the equivalence $\abs{x-z} <\abs{x+z} \iff z\cdot x>0$. 
Since  $\bar \psi_R$ and $\bar \vfi_R$ are non-increasing, we deduce that $\grad A(x) \cdot x \leq0$, this is,  $A$ is radially non-increasing. 

The general case $\psi\in L^1(\R^N)$ is then proven by approximating $\psi$ by $\psi_n\in W^{1,1}(\R^N)$, which satisfy the property $A_n(x)\leq A_n(y)$ if $\abs{x}\geq\abs{y}$. Then, the pointwise convergence of a subsequence $\psi_{n_j}\to \psi$ implies that this property is kept by the limit, $A$. $\hfill\Box$

\section{Proof of Theorem~\ref{th:evolution}}\label{sec:evolution}

Consider the partition of the interval $[0,T]$ given by 
$t_n=n\tau$, for $n=0,\ldots,N$, with $\tau = T/N$.
We introduce the following explicit time discretization of the problem ENP($E,h,\zeta,\eps$). We set 
$w_0=\zeta$ in $E$ and, for $n=0,\ldots,N-2$,
\begin{align*}
w_{n+1}(x)= (1-\tau c)w_n(x) + \tau\int_{\R^N} J_\eps(x-y)(w_n(y)-w_n(x)) dy +\tau h_n(x) 
\end{align*}
for $x\in E$, and $w_{n+1}(x)=0$ for $x\in \R^N\backslash E$, where $h_n(x)=h(t_n,x)$.
We also consider the piecewise constant interpolant
\begin{align*}
 w^{(\tau)}(t,x)=w_{n}(x)\qtextq{for} (t,x)\in (t_n,t_{n+1}]\times E \qtextq{and}n=0,\ldots,N-1.
\end{align*}
Let $v$ and $u$ be the solutions of ENP($\O,G,v_0,\eps$) and ENP($\O^*,G^*,v_0^*,\eps$), and let $v^{(\tau)}$ and $u^{(\tau)}$ be their corresponding approximations according to the discrete scheme. Here and in what follows,  $G^*$ denotes the Schwarz symmetrization of $G$ with respect to the space variable, i.e. $(G(t,\cdot))^*$.

\medskip
\no \emph{Assumption (H)$_4$:
The data, $\O\subset\R^N$, $G:(0,T)\times\O\to\R$, and $v_0:\O\to\R$ are regular enough to provide the convergence 
\begin{align}
&v^{(\tau)} \to v \qtext{strongly in }L^p((0,T)\times \O) \text{ and a.e. in }(0,T)\times \O, \label{h4.1}\\
&u^{(\tau)} \to u \qtext{strongly in }L^p((0,T)\times \O^*) \text{ and a.e. in }(0,T)\times \O^*.\label{h4.2}
\end{align} 
}

\begin{remark}
	The convergence property assumed in (H)$_4$ has been proven for a general class of nonlinear nonlocal diffusion evolution problems. See \cite{Galiano2019, Galiano2021} for details.
\end{remark}

The semi-discrete problem may be reformulated as: $w_0=\zeta$ in $E$ and, for $n=0,\ldots,N-2$, 
\begin{align}
\label{esc.ex2}
w_{n+1}(x) = & \alpha \int_{\R^N} \rho_\eps (x-y) w_{n}(y)  dy + \beta w_n(x)+ \tau h_n(x)\qtextq{for}x\in E,
\end{align}
and $w_{n+1}=0$ for $x\in \R^N\backslash E$, where
\begin{align*}
	\alpha = \frac{\tau}{\eps^2},\quad  \rho_\eps  = \eps^2J_\eps,\quad \beta=(1-\tau(c+\eps^{-2})).
\end{align*}
We choose the time step $\tau<1/(c+\eps^{-2})$, implying $\beta>0$. Observe that this is the usual stability condition associated to explicit time discretizations.

Like in the proof of Theorem~\ref{th:stationary}, we remove from the notation the reference to the parameter $\eps$ and refer to  the previous scheme as to  DNS($E,h,\zeta$).

%

Let $v_n$ and $u_n$ be the sequences defined by the schemes DNS($\O,G,v_0$) and DNS($\O^*,G^*,v_0^*$), respectively. From  \fer{esc.ex2}, it is clear that, as an addition of non-negative radially non-increasing functions, $u_n$ is radially non-increasing.

Multiplying the identity \fer{esc.ex2} corresponding to  DNS($\O,G,v_0$) by any non-negative  $\vfi\in L^{p'}(\O)$, 
integrating in $\O$, and using the symmetry of $\rho $ to interchange the convolution of $\rho $ with $v_n$ for the convolution of $\rho $ with $\vfi$, and taking into account that $v_n=0$ in $\R^N\backslash\O$, we get 
\begin{align*}
&\int_{\O} v_{n+1} \vfi = \alpha \int_{\O} \int_{\O} \rho (x-y) \vfi(y) dy  v_n(x)  dx + \int_{\O}  (\beta v_n + \tau G_n)\vfi.
\end{align*}
Proceeding similarly for DNS($\O^*,G^*,v_0^*$), we obtain
\begin{align*}
&\int_{\O^*} u_{n+1} \vfi^* = \alpha \int_{\O^*} \int_{\O^*} \rho (x-y) \vfi^*(y) dy  u_n(x)  dx + \int_{\O^*}   (\beta u_n + \tau G^*_n)\vfi^*.
\end{align*}
For $n=0$, we have
\begin{align}
&\int_{\O} v_{1} \vfi = \alpha \int_{\O} \int_{\O} \rho (x-y) \vfi(y) dy  v_0(x)  dx + \int_{\O}  (\beta v_0 + \tau G_0)\vfi. \label{ncero}
\end{align}
Riesz's inequality and $u_0=v_0^*$ imply
\begin{align*}
 \int_{\O} \int_{\O} \rho (x-y)  \vfi(y) dy  v_0(x)  dx \leq  \int_{\O^*} \int_{\O^*} \rho (x-y) \vfi^*(y) dy  u_0(x)  dx.
\end{align*}
Hardy-Littlewood's inequality gives 
\begin{align*}
\int_{\O}  v_0\vfi \leq \int_{\O^*}  u_0\vfi^* \qtextq{and}
\int_{\O}  G_0\vfi \leq \int_{\O^*}  G_0^*\vfi^*.
\end{align*}
Therefore, \fer{ncero} may be estimated as 
\begin{align}
\int_{\O} v_{1} \vfi &\leq \alpha \int_{\O^*} \int_{\O^*} \rho (x-y) \vfi^*(y) dy  u_0(x)  dx +  \int_{\O^*}  (\beta u_0+\tau G_0^*)\vfi^* \nonumber \\
&= \int_{\O^*} u_{1} \vfi^* \label{est.r1}
\end{align}
Similarly, for $n=1$, we have 
\begin{align}
\label{est.r2}
\int_{\O} v_{2} \vfi &= \alpha \int_{\O} \int_{\O} \rho (x-y) \vfi(y) dy  v_1(x)  dx +\int_{\O}  (\beta v_1 + \tau G_1)\vfi.
\end{align}
Hardy-Littlewood's inequality and  \fer{est.r1} give
\begin{align}
\label{est.r3}
	\int_{\O}  (\beta v_1 + \tau F_1)\vfi \leq \int_{\O^*}  (\beta u_1 + \tau G_1^*)\vfi^*.
\end{align}
Replacing the test function $\vfi$ by $\rho \co_\O \vfi$ in \fer{ncero}, we get 
\begin{align}
\int_{\O} v_{1}(x)&\int_{\O} \rho (x-y) \vfi(y) dydx \nonumber \\
&= \alpha \int_{\O} \int_{\O} \rho (x-y)  \Big(\int_{\O} \rho (y-z) \vfi(z) dz \Big) dy  v_0(x)  dx \nonumber\\&
+ \int_{\O}  (\beta v_0 + \tau G_0) \Big(\int_{\O} \rho (x-y) \vfi(y) dy \Big)dx.\label{est.r6}
\end{align}
The first term of the right hand side is estimated by means of the generalized Riesz's inequality, see \cite[Theorem~8]{Lieb2001}, as 
\begin{align*}
 \int_{\O}\int_{\O}\int_{\O} v_0(x) \vfi(z) &\rho (x-y) \rho (y-z) dx dy dz  \\
& \leq  \int_{\O^*}\int_{\O^*}\int_{\O^*} u_0(x) \vfi^*(z) \rho (x-y) \rho (y-z) dx dy dz ,
\end{align*}
while for the second we use Riesz's inequality to obtain
\begin{align*}
	\int_{\O}  v_0(x) \int_{\O} \rho (x-y) \vfi(y) dy dx &\leq
	\int_{\O^*}  u_0(x) \int_{\O^*} \rho (x-y) \vfi^*(y) dy dx, \\
	\int_{\O}  G_0(x) \int_{\O} \rho (x-y) \vfi(y) dy dx &\leq
	\int_{\O^*}  G_0^*(x) \int_{\O^*} \rho (x-y) \vfi^*(y) dy dx. 
\end{align*}
Thus, \fer{est.r6} yields 
\begin{align}
\label{est.r5}
\int_{\O} v_{1}(x) \Big(\int_{\O} \rho (x-y) \vfi(y) dy \Big)dx \leq 
\int_{\O^*} u_{1}(x) \Big(\int_{\O^*} \rho (x-y) \vfi^*(y) dy \Big)dx.
\end{align}
Returning to \fer{est.r2} and taking into account  \fer{est.r3} and  \fer{est.r5}, we obtain 
\begin{align*}
\int_{\O} v_{2}(x) \vfi(x)dx \leq \int_{\O^*} u_{2}(x) \vfi^*(x)dx.
\end{align*}
It is now clear that we may repeat this argument recurrently to deduce 
\begin{align*}
\int_{\O} v_{n} \vfi \leq \int_{\O^*} u_n \vfi^*, \qtextq{for all}n=0,1,\ldots,N-2.
\end{align*}
By choosing $\vfi=1$, for $p=1$ and $\vfi=v_n^{p-1}$, for $p>1$, we deduce like in \fer{p1:7} that  $\norm{v_n}_{L^p(\O)} \leq \norm{u_n}_{L^p(\O^*)}$, for all $ n=0,1,\ldots,N-2$, implying that  $\norm{v^{(\tau)}}_{L^p((0,T)\times \O)} \leq \norm{u^{(\tau)}}_{L^p((0,T)\times\O^*)}$.

Finally, using the assumption (H)$_4$ and the argument of the proof of the corollaries, see Remark~\ref{rem:cor}, we deduce \fer{th2:comp}.

$\hfill\Box$

\end{document}